\documentclass{article}

\usepackage[utf8]{inputenc}

\usepackage{amsmath}
\usepackage{amssymb}
\usepackage{vmargin}

\usepackage{parskip}
\usepackage{subfig}
\usepackage{graphicx}

\usepackage{cite}

\usepackage{authblk}

\usepackage{color}
\title{Number-conserving cellular automata with a von Neumann neighborhood of range one}
\author{Barbara Wolnik}
\author{Adam Dzedzej}
\affil{Institute of Mathematics, Faculty of Mathematics, Physics and Informatics, University of Gdańsk, 80-308 Gdańsk, Poland} 
\author{Jan M.~Baetens}
\author{Bernard De Baets}
\affil{KERMIT, Department of Mathematical Modelling, Statistics and Bioinformatics, Ghent University, Ghent, Belgium}

\date{}

\newcommand{\ff}[5]{
\setlength{\arraycolsep}{1pt}
\renewcommand{\arraystretch}{0.8}
{f}\!\left(\!
\begin{array}{rcl}
 &#1&  \\
 #2 & #3 & #4 \\
 & #5 & 
\end{array}\!
\right)
}

\newcommand{\fff}[7]{
\setlength{\arraycolsep}{1pt}
\renewcommand{\arraystretch}{0.8}
{f}\!\left(\!
\begin{array}{rcl}
 &#1& {\color{blue}#7} \\
 #2 & #3 & #4 \\
{\color{blue}#6} & #5 &  
\end{array}\!
\right)
}

\newcommand{\seff}[5]{
\setlength{\arraycolsep}{1pt}
\renewcommand{\arraystretch}{0.8}
f^{E}\!\left(\!
\begin{array}{rcl}
 &#1&  \\
 #2 & #3 & #4 \\
 & #5 & 
\end{array}\!
\right)
}

\newcommand{\dimer}[2]{
D_{\substack{#1\\#2}}
}

\newcommand{\vv}[2]{
{ \overrightarrow{{\bf #1}_{#2}}}
}

\newcommand{\C}{\mathcal{C}}
\newcommand{\NN}{\mathcal{N}}

\newcommand{\bx}{\mathbf{x}}
\newcommand{\bze}{\mathbf{0}}
\newcommand{\bi}{\mathbf{i}}
\newcommand{\bj}{\mathbf{j}}
\newcommand{\bl}{\mathbf{l}}
\newcommand{\bq}{\mathbf{q}}

\newcommand{\Z}{\mathbb{Z}}

\newcommand{\R}{\mathbb{R}}
\newcommand{\PP}{\mathbf{\Omega}}
\newcommand{\CP}{\mathbf{\Lambda}}
\newcommand{\di}{\operatorname{dist}}

\newtheorem{thm}{Theorem}[section]
\newtheorem{defn}[thm]{Definition}
\newtheorem{lem}[thm]{Lemma}
\newtheorem{cor}[thm]{Corollary}

\newtheorem{rem}[thm]{Remark}

\usepackage{letltxmacro}
\LetLtxMacro{\originaleqref}{\eqref}
\renewcommand{\eqref}{Eq.\ \originaleqref}

\begin{document}

\maketitle

\begin{abstract}
 We present necessary and sufficient conditions for a cellular automaton with a von Neumann neighborhood of range one to be number-conserving. The conditions are formulated for any dimension  and for any set of states containing zero. The use of the geometric
 structure of the von Neumann neighborhood allows for computationally
 tractable conditions even in higher dimensions. 
\end{abstract}

\section{Introduction}
This paper focuses on $d$-dimensional cellular automata (CAs) that possess the most popular additive invariant, namely the sum of the states of all cells. This kind of CA, called number-conserving, has been introduced by Nagel and Schreckenberg~\cite{NS} in the early nineties and has received ample attention in the literature, especially as a model of a system of interacting particles moving in a lattice. In particular, such CAs occur naturally in the context of highway traffic~\cite{Belitsky2001, KKW02, PhysRevLett.90.088701} and fluid flow~\cite{PhysRevLett.56.1505}. In one dimension, necessary and sufficient conditions 
for a CA to be number-conserving are given by Boccara and Fukś in~\cite{BF98} for two states ($Q=\{0,1\}$) and in~\cite{BoccaraF02} for $q$ states ($Q=\{0,1,\ldots ,q-1\}$). The formulas in these works can be also derived from the Hattori and Takesue theorem~\cite{HT} concerning general additive conserved quantities. The approach of Boccara and Fukś allows to enumerate all number-conserving one-dimensional CAs in the case of few states and a small neighborhood.

Durand et al.~\cite{Durand2003} state necessary and sufficient conditions for a CA to be number-conserving in two or more dimensions for $q$ states ($Q=\{0,1,\ldots,q-1\}$). 
They also formalized three notions of number conservation (periodic, finite, infinite) and proved their equivalence. These results were obtained for product neighborhoods, so it is universal and shows, for example, that number conservation is a decidable property (provided that the number of states is finite). Moreira~\cite{Mo3} extended the work of Durand et al.\ to any finite subset of $\Z$ as set of states. Unfortunately, the formulas in~\cite{Durand2003} %
% * <bernard.debaets@ugent.be> 2017-04-29T14:28:34.745Z:
%
% >  product neighborhoods, so it is universal and shows, for example,
%
% ^.
and~\cite{Mo3} are too complicated to be effective in enumerating all number-conserving two- or higher-dimensional CAs, even in the case of two states. Furthermore, in these approaches a CA with the von Neumann neighborhood was considered as a CA with the Moore neighborhood, 
which implies an unnecessary complication and strongly increases the computational complexity.

Tanimoto and Imai~\cite{TI} provided a characterization of two-dimensional number-conserving CAs with the von Neumann neighborhood. Their result is stated in terms of some flow functions (in the vertical, horizontal and diagonal direction), and allows, for example, to prove facts about rotation-symmetric CAs. They showed that rotation-symmetric number-conserving CAs with at most four states are trivial and also gave a full characterization of five-state rotation-symmetric number-conserving CAs.
% * <bernard.debaets@ugent.be> 2017-04-29T14:31:15.064Z:
%
% > rotation-symmetric CAs.
%
% ^.

In this paper, we restrict our attention to the simplest but most commonly used neighborhood: the von Neumann neighborhood of range one, referred to as ``the von Neumann neighborhood" for the sake of simplicity. We formulate necessary and sufficient conditions for $d$-dimensional CAs with the von Neumann neighborhood to be number-conserving in a way similar as in~\cite{BoccaraF02}. The conditions apply for any state set (finite or not), but for our convenience, we assume that the set of states includes $0$.
In fact, the only property of $0$ we use is its quiescence, so $0$ can be replaced by any other state, since all states of a number-conserving CA are quiescent (see~Lemma \ref{lemma:qq}). 
% * <bernard.debaets@ugent.be> 2017-04-29T14:33:53.985Z:
%
% > In fact, the only property of $0$ we use is its quiescence, so $0$ can be replaced by any other state, since all states of a number-conserving CA are quiescent (see~Lemma \ref{lemma:qq}). 
%
% ^.
The form of our characterization allows to enumerate all number-conserving $d$-dimensional CAs in the case of few states and small~$d$.

This paper is organized as follows. In Section 2 the basic concepts and notations are introduced. Section 3 presents the necessary and sufficient conditions
for a two-dimensional CA to be number-conserving, while the higher-dimensional case is considered in Section 4. Presenting the results in this way, guarantees that they are 
accessible for any reader. The final section summarizes the main results and gives some hints for further research. 

\section{Preliminaries}

In this section, we introduce CAs and recall some results 
that we will use in the following sections. To define a~CA, 
one needs to specify a~space of cells, a neighborhood and a local rule.

\subsection{The cellular space}\label{cs}

Let us fix the dimension $d\geq 1$ and consider natural numbers $n_1,n_2,\ldots,n_d$ greater than $4$. We consider the cellular space as a grid with periodic boundary conditions, defined by
\[
\C=\left(\Z / n_1\Z\right)\times \left(\Z / n_2\Z\right)\times\ldots\times \left(\Z / n_d\Z\right)
\]
\[
=\{0,1,\ldots,n_1-1\}\times\{0,1,\ldots,n_2-1\}\times\ldots \times\{0,1,\ldots,n_d-1\} .
\]
With this notation, each cell $\bi\in\C$ is a $d$-tuple 
$(i_1,\ldots, i_d)$, where $i_k\in\Z / n_k\Z$. In the case $d=2$, we prefer to denote cells by $(i,j)$ rather than $(i_1, i_2)$. 
Denoting the number of elements of a set $A$ as $|A|$, we have $|\C|=n_1\cdot n_2\cdot\ldots \cdot n_d$. 

\begin{figure}[!ht]
\begin{center}
 \subfloat[]{
\includegraphics[height=5cm]{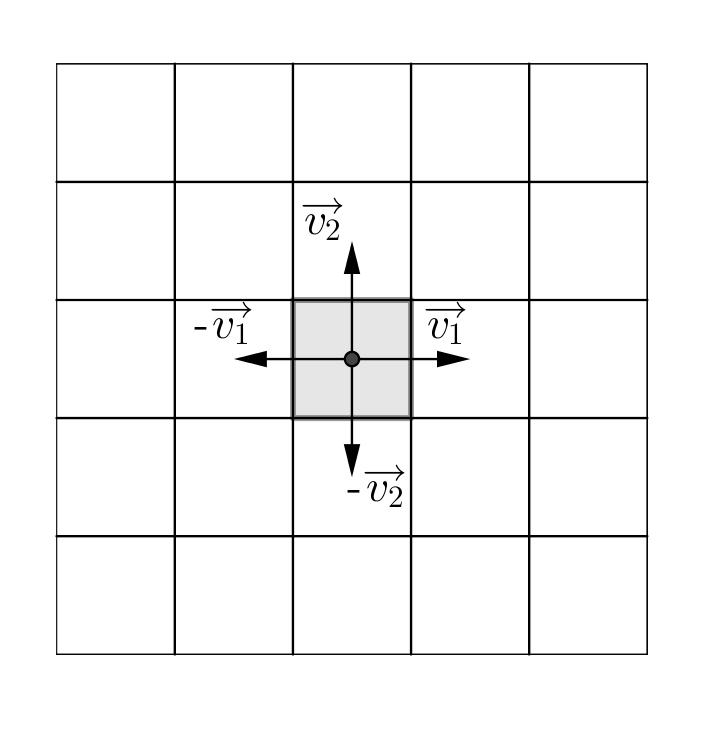}
}
 \subfloat[]{
\includegraphics[height=50mm]{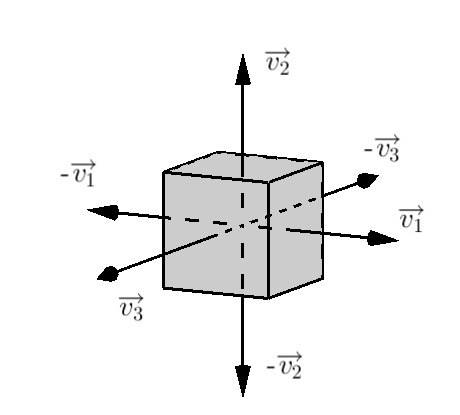}
}

\end{center}
\caption{Neighborhood directions in the case of two- (a) and three- (b) dimensional CAs.}
\label{neighborhoodFig}
\end{figure}

Due to the periodic boundary conditions, each cell in $\C$ has exactly $2d$ adjacent cells: two in each of 
the orthogonal axis directions. For example, if $d=2$, then there are four adjacent cells: two in the horizontal direction, 
$\vv{v}{1}$ (right) and $\mbox{-}\vv{v}{1}$ (left), and two in the vertical direction, $\vv{v}{2}$ (up) and $\mbox{-}\vv{v}{2}$ (down), as shown
in Fig.~\ref{neighborhoodFig}(a).
If $d=3$, then there are two additional adjacent cells: in the directions $\vv{v}{3}$ (forward) and $\mbox{-}\vv{v}{3}$ (backward),
as shown in Fig.~\ref{neighborhoodFig}(b).

To be able to describe the situation in general, \emph{i.e.}\ for any dimension $d$, we introduce the following notation.
For each $k\in\{1,2,\ldots,d\}$, we define
the vector  $\vv{v}{k}=(0,0,\ldots,0,1,0,\ldots,0)\in \R^d$, where the $k$-th component is equal to $1$ and all others are equal to zero. Let us denote the set of all considered directions as $V_+$, \emph{i.e.}\
\[
V_+=\{\vv{v}{1},\mbox{-}\vv{v}{1},\vv{v}{2},\mbox{-}\vv{v}{2},\ldots, \vv{v}{d},\mbox{-}\vv{v}{d}\}.
\]
Additionally, let $\vv{0}{}=(0,0,\ldots,0)\in \R^d$ and
$V= V_+\cup \{\vv{0}{} \}$.

 For $\bi \in \C$ and $\vv{v}{}\in V$, we now define the sum $\bi +\vv{v}{}\in \C$ as a cell adjacent to $\bi$ in direction $\vv{v}{}$, if
 $\vv{v}{}\in V_+$ or as $\bi$ itself if $\vv{v}{}=\vv{0}{}$.

\subsection{The neighborhood}\label{sect:neibor}

As we mentioned in Section 1,  we only consider the von Neumann neighborhood. For each cell $\bi \in\C$, its von Neumann neighborhood $P(\bi)$  consists of the cell $\bi$ and 
its $2d$ adjacent cells:
\[
P(\bi )= \bi + V := \{ \bi + \vv{v}{}\mid \vv{v}{}\in V \}\,.
\]
This neighborhood can be described using the Manhattan distance defined by
\begin{equation}\label{di}
\di (\bi ,\bj )=\sum_{k=1}^{d}\min\left( |i_k-j_k|, n_k-|i_k-j_k|\right)\,,
\end{equation}
where $\bi =(i_1,i_2,\ldots,i_d)\in\C$ and $\bj =(j_1,j_2,\ldots,j_d)\in\C$.
From this point of view, the von Neumann neighborhood $P(\bi)$ consists of those cells whose distance from $\bi$ is not greater than $1$, \emph{i.e.}\
\[
P(\bi )=\{ \bj \in\C\mid \di (\bi ,\bj )\leq 1\},. 
\]
We now describe some relations between $P(\bi)$ and $P(\bj)$ in case $\bi\neq \bj$.

First of all, we observe that if $\di (\bi ,\bj )>2$, then $P(\bi )\cap P(\bj )=\emptyset$. Indeed, if there is some cell $\bl \in P(\bi )\cap P(\bj )$, then 
$\di (\bi ,\bl )\leq 1$ and  $\di (\bj ,\bl )\leq 1$, which implies that $\di (\bi ,\bj )\leq \di (\bi ,\bl )+\di (\bl ,\bj )\leq 2$. 
Secondly, if $\di (\bi ,\bj )=1$, then cells $\bi$ and $\bj$ are adjacent. 
Hence, $\bj \in P(\bi )$ and $\bi \in P(\bj )$ and there are no other cells in $P(\bi )\cap P(\bj )$, which implies that  $\bj =\bi +\vv{v}{}$ for some $\vv{v}{} \in V_+$. Finally, the case $\di (\bi ,\bj )=2$ is the most interesting one. The size of $P(\bi )\cap P(\bj )$ depends on the relative positions of $\bi$ and $\bj$, as we can either move from $\bi$ to $\bj$ in two steps in the same direction or one step in one direction and a second step in some orthogonal direction.
Hence, there are two cases:
\begin{itemize}
\item[1)] $\bj = \bi + \vv{v}{} + \vv{v}{}$, for some $\vv{v}{}\in V_+$, so  $P(\bi )\cap P(\bj )=\{\bi +\vv{v}{}\}$.
\item[2)] $\bj = \bi +\vv{u}{} + \vv{v}{}$, for some $\vv{u}{},\vv{v}{}\in V_+$, where both $\vv{u}{}\neq \vv{v}{}$ and $\vv{u}{}\neq \mbox{-}\vv{v}{}$. In this case, $P(\bi )\cap P(\bj )=\{\bi+\vv{u}{}, \bi +\vv{v}{}\}$. 
\end{itemize}

The following lemma summarizes these observations.
\begin{lem}\label{lemma:vN}
Let $\bi ,\bj \in\C$ and $\bi \neq \bj $. 
\begin{itemize}
\item[(a)] If $\bj =\bi +\vv{v}{}$ for some $\vv{v}{}\in V_+$, then $P(\bi )\cap P(\bj )=\{ \bi+\vv{0}{} ,\bi+\vv{v}{}\}$.
\item[(b)] If $\bj = \bi +\vv{v}{} + \vv{v}{}$ for some $\vv{v}{}\in V_+$, then $P(\bi )\cap P(\bj )=\{ \bi+\vv{v}{} \}$. 
\item[(c)] If $\bj = \bi +\vv{u}{} + \vv{v}{}$
 for some $\vv{u}{},\vv{v}{}\in V_+$ where both $\vv{u}{}\neq \vv{v}{}$ and 
 $\vv{u}{}\neq \mbox{-}\vv{v}{}$, 
 then $P(\bi )\cap P(\bj)=\{\bi+\vv{u}{}, \bi +\vv{v}{}\}$.
\item[(d)] In all other cases, it holds that $P(\bi )\cap P(\bj )=\emptyset$.
\end{itemize}
\end{lem}

We now define the set $\PP$, containing all possible pairs of vectors from $V_+$ used in Lemma~ \ref{lemma:vN} in the descriptions of $P(\bi )\cap P(\bj)$ in cases (a) and (c). Thus $\PP$ consists of pairs
 $\lbrace \vv{0}{}, \vv{v}{} \rbrace$,
 for every $\vv{v}{}\in V_+$, and pairs 
 $\lbrace \vv{u}{}, \vv{v}{}\rbrace$, for every $\vv{u}{},\vv{v}{}\in V_+$ such that $\vv{u}{}\neq \vv{v}{}$ and $\vv{u}{}\neq \mbox{-}\vv{v}{}$. 
 
As the pairs  $\lbrace \vv{u}{}, \vv{v}{}\rbrace$ and $\lbrace \vv{v}{}, \vv{u}{}\rbrace$ are equal, the set $\PP$ contains exactly $2d^2$ elements. 
For example, if $d=2$, the set $\PP$ contains eight elements:
\begin{equation}
\PP = \bigg\lbrace  
\lbrace \vv{0}{},\vv{v}{1}\rbrace,  \lbrace\vv{0}{},\vv{v}{2}\rbrace ,  \lbrace \vv{v}{1},\vv{v}{2}\rbrace, 
\lbrace\vv{v}{1},\mbox{-}\vv{v}{2}\rbrace, 
\lbrace \vv{0}{},\mbox{-}\vv{v}{1}\rbrace,  \lbrace\vv{0}{},\mbox{-}\vv{v}{2}\rbrace,  \lbrace\mbox{-}\vv{v}{1},\mbox{-}\vv{v}{2}\rbrace, \lbrace\mbox{-}\vv{v}{1},\vv{v}{2}\rbrace \bigg\rbrace\,,
\label{d2vectors}
\end{equation}
while if $d=3$, it holds that $|\PP|=18$:
\begin{equation}
\setlength{\arraycolsep}{0pt}
\begin{array}{l}
\PP = \bigg\lbrace
\lbrace\vv{0}{},\vv{v}{1}\rbrace,  \lbrace\vv{0}{},\vv{v}{2}\rbrace, 
 \lbrace\vv{0}{},\vv{v}{3}\rbrace,  \lbrace\vv{v}{1},\vv{v}{2}\rbrace,  \lbrace\vv{v}{1},\mbox{-}\vv{v}{2}\rbrace,  \lbrace\vv{v}{1},\vv{v}{3}\rbrace,  \lbrace\vv{v}{1},\mbox{-}\vv{v}{3}\rbrace,  \lbrace\vv{v}{2},\vv{v}{3}\rbrace,  \lbrace\vv{v}{2},\mbox{-}\vv{v}{3}\rbrace,\\
\lbrace\vv{0}{},\mbox{-}\vv{v}{1}\rbrace,  \lbrace\vv{0}{},\mbox{-}\vv{v}{2}\rbrace,  \lbrace\vv{0}{},\mbox{-}\vv{v}{3}\rbrace,  \lbrace\mbox{-}\vv{v}{1},\mbox{-}\vv{v}{2}\rbrace, \lbrace\mbox{-}\vv{v}{1},\vv{v}{2}\rbrace,  \lbrace\mbox{-}\vv{v}{1},\mbox{-}\vv{v}{3}\rbrace,  \lbrace\mbox{-}\vv{v}{1},\vv{v}{3}\rbrace,  \lbrace\mbox{-}\vv{v}{2},\mbox{-}\vv{v}{3}\rbrace,  \lbrace\mbox{-}\vv{v}{2},\vv{v}{3}\rbrace
\bigg\rbrace\,.
\end{array}
\label{d3vectors}
\end{equation}
 From Lemma~\ref{lemma:vN}, we deduce the following.
\begin{rem}\label{remark}
Let $\bi\in\C$. If for some $\bj\in\C$ it holds that $|P(\bi ) \cap P(\bj )| =2$,
then there exists a unique element $\lbrace\vv{u}{}, \vv{v}{}\rbrace\in\PP$ such that 
$\bj =(\bi +\vv{u}{})+\vv{v}{}$. Moreover, one of the shared cells is $\bi +\vv{u}{} = \bj +(\mbox{-}\vv{v}{})$, while the other one is $\bi +\vv{v}{}
= \bj + (\mbox{-}\vv{u}{})$.
\end{rem}
For a given pair $\lbrace \vv{u}{}, \vv{v}{}\rbrace\in\PP$, the pair $\lbrace\mbox{-}\vv{u}{},
 \mbox{-}\vv{v}{}\rbrace$ will be called its~\emph{matching pair}. For example, in the two-dimensional case,
$\lbrace\vv{0}{},\vv{v}{1}\rbrace, \lbrace\vv{0}{},\vv{v}{2}\rbrace, \lbrace\vv{v}{1},\vv{v}{2}\rbrace, \lbrace\vv{v}{1},\mbox{-}\vv{v}{2}\rbrace$ are the matching pairs of
$\lbrace \vv{0}{},\mbox{-}\vv{v}{1}\rbrace,  \lbrace\vv{0}{},\mbox{-}\vv{v}{2}\rbrace, \lbrace\mbox{-}\vv{v}{1},\mbox{-}\vv{v}{2}\rbrace,  \lbrace\mbox{-}\vv{v}{1},\vv{v}{2}\rbrace$, 
respectively, and vice versa.

Obviously, different pairs in $\PP$ have different matching pairs. If from each of the $d^2$ pairs of matching pairs we select one pair, then we obtain a set denoted by $\CP$.
We can construct $\CP$ in $2^{d^2}$ ways, but it always holds that 
 $|\CP |=d^2$.

\subsection{The configuration space}

Here, we consider an arbitrary non-singleton set $Q\subseteq\R$ containing zero and define 
$Q_+:= Q\setminus\{0\}$.
By a configuration, we mean any mapping from the grid $\C$ to $Q$. The set of all possible configurations is denoted by $X=Q^{\C}$. 
The state of cell $\bi $ in a configuration $\bx\in X$ is denoted by $\bx(\bi )$ or, if $d=2$, by $x_{i,j}$ for the cell $(i,j)$. 

Given a configuration $\bx\in X$, we define the sum of the states in $\bx$ as:
\[\sigma(\bx)=\sum_{\bi \in\C}\bx(\bi )\,.
\]

\subsection{The set of neighborhood configurations}

By the set of all possible neighborhood configurations $\NN$, we mean the set of all functions $N:\; V\to Q$. If $N$ is identically equal to zero, then we call it \emph{trivial}. As the set $V$ has $2d+1$ elements, namely $\vv{0}{}$, $\vv{v}{1}$, $\mbox{-}\vv{v}{1}$, $\ldots$, $\vv{v}{d}$, $\mbox{-}\vv{v}{d}$, we can define any neighborhood configuration by the sequence $(N(\vv{0}{}), N(\vv{v}{1}), N(\mbox{-}\vv{v}{1}), \ldots, N(\vv{v}{d}), N(\mbox{-}\vv{v}{d}))$. In the two-dimensional case, for mnemotechnical reasons, we represent $N$ as 
$
\setlength{\arraycolsep}{1pt}
\renewcommand{\arraystretch}{0.8}
\begin{array}{rcl}
 & q_1 &  \\
q_2 & q_3 & q_4, \\
 & q_5 & 
\end{array}$
with $N(\vv{0}{})=q_3$, $N(\vv{v}{1})=q_4$, $N(\mbox{-}\vv{v}{1})=q_2$, $N(\vv{v}{2})=q_1$ and $N(\mbox{-}\vv{v}{2})=q_5$. In this way,
the set $\NN$ can be defined as
\[
\NN=\left\{ 
\setlength{\arraycolsep}{1pt}
\renewcommand{\arraystretch}{0.8}
\begin{array}{rcl}
 & q_1 &  \\
q_2 & q_3 & q_4 \\
 & q_5 & 
\end{array}\mid q_1,q_2,q_3,q_4,q_5\in Q
\right\}.
\]

For $d=3$, we similarly have $\setlength{\arraycolsep}{1pt}
\renewcommand{\arraystretch}{0.8}
\begin{array}{rcl}
 & q_1 & {\color{blue}q_7} \\
 q_2 & q_3 & q_4, \\
{\color{blue}q_6} & q_5 &  
\end{array}$ with $N(\vv{0}{})=q_3$, $N(\vv{v}{1})=q_4$, $N(\mbox{-}\vv{v}{1})=q_2$, $N(\vv{v}{2})=q_1$, $N(\mbox{-}\vv{v}{2})=q_5$, $N(\vv{v}{3})=q_6$ and $N(\mbox{-}\vv{v}{3})=q_7$. 

Let $q\in Q$. By $H_{q}$ we denote the homogeneous neighborhood configuration, \emph{i.e.}
\[
(\forall \vv{v}{}\in V) (H_{q}(\vv{v}{})=q)\,.
\]
Let $\vv{v}{}\in V$ and $q\in Q$ be given. By $M_{\vv{v}{}:q}$ we denote the neighborhood configuration which differs from $H_{0}$ in at most one component. More precisely, $M_{\vv{v}{}:q}$ takes the value $q$ in direction $\vv{v}{}$ and zero in the other directions, \emph{i.e}
\[
\left( \forall \vv{u}{}\in V \right) \left(M_{\vv{v}{}:q}(\vv{u}{})=
\begin{cases}
q \quad\text{, if }  \vv{u}{}= \vv{v}{}\\
0 \quad\text{, if }  \vv{u}{}\neq \vv{v}{}
\end{cases}\right)\,.
\]
Every $M_{\vv{v}{}:q}$ is called a \emph{monomer}.
Of course, if $q=0$, then $M_{\vv{v}{}:q}$ trivially equals $H_0$. 

Similarly, if $(\vv{u}{},\vv{w}{})\in \PP$ and $p,q\in Q$, then $\dimer{\vv{u}{}:p}{\vv{w}{}:q}$ denotes the neighborhood configuration that takes the value $p$ in direction $\vv{u}{}$, the value $q$ in direction $\vv{w}{}$ and zero in the other directions, \emph{i.e.}
\[
\left(
\forall \vv{v}{}\in V \right) \left(\dimer{\vv{u}{}:p}{
\vv{w}{}:q}(\vv{v}{})=
\begin{cases}
p \quad\text{, if }\vv{v}{}=\vv{u}{}\\
q \quad\text{, if }\vv{v}{}=\vv{w}{}\\
0 \quad\text{, otherwise}
\end{cases} \right).
\]
It is obvious that $\dimer{\vv{u}{}:p}{\vv{w}{}:q}$ equals $\dimer{\vv{w}{}:q}{\vv{u}{}:p}$.
A neigborhood of the type $\dimer{\vv{u}{}:p}{
\vv{w}{}:q}$ is called a {\it dimer}. As the pairs $\lbrace\vv{u}{}, \vv{w}{}\rbrace$ and $\lbrace \mbox{-}\vv{u}{}, \mbox{-}\vv{w}{}\rbrace$ 
are matching, we also refer to the dimers $\dimer{\vv{u}{}:p}{\vv{w}{}:q}$ and $\dimer{\mbox{-}\vv{w}{}:p}{\mbox{-}\vv{u}{}:q}$ as matching dimers.
%We call these matching, since if $\dimer{\vv{u}{}:p}{\vv{w}{}:q}$ is given, then by Lemma \ref{lemma:eq} or Lemma \ref{lemma:cp3} we can calculate the value of $\dimer{\mbox{-}\vv{w}{}:p}{\mbox{-}\vv{u}{}:q}$ using only monomers.

\noindent
If $\bx\in X$ and $\bi \in\C$ are given, then $N_{\bx,\bi}$ denotes the configuration of the von Neumann neighborhood of cell $\bi$ in the configuration $\bx$, \emph{i.e.}
\[
\left(\forall \vv{v}{}\in V \right) \left( N_{\bx,\bi} (\vv{v}{})=\bx(\bi +\vv{v}{})\right)\,.
\]

\subsection{Local and global rules}
A function $f\colon \NN\to Q$ is called a {\it local} rule.  The {\it monomer expansion} $f^{E}\colon \NN\to\R$ of a local rule $f$ is defined by
\[
f^{E}(N)=\sum_{\vv{v}{}\in V}f\left( M_{\vv{v}{}:N(\vv{v}{})}\right)\,.
\]
Note that $M_{\vv{v}{}:N(\vv{v}{})}$ is obtained from $N$ by keeping the value in the direction $\vv{v}{}$, while setting all other values to zero.
Thus, in the monomer expansion, the neighborhood configuration is decomposed into monomers, the local rule $f$ is applied to each monomer separately the resulting valued are added. For example, if $d=2$, then we have
\[
f^{E}\left(
\setlength{\arraycolsep}{1pt}
\renewcommand{\arraystretch}{0.8}
\begin{array}{rcl}
 & q_1 &  \\
q_2 & q_3 & q_4 \\
 & q_5 & 
\end{array}
\right)
= \ff{q_1}{0}{0}{0}{0}+\ff{0}{q_2}{0}{0}{0}+\ff{0}{0}{q_3}{0}{0}+\ff{0}{0}{0}{q_4}{0}+\ff{0}{0}{0}{0}{q_5}\,.
\]
Similarly, in the $d$-dimensional case, we have
\begin{equation}                                                                             \label{dimerexp}
f^{E}\left(\dimer{\vv{u}{}:p}{\vv{w}{}:q}\right) = 
f\left( M_{\vv{u}{}:p}\right)
+f\left( M_{\vv{w}{}:q}\right)
\end{equation}
for any $(\vv{u}{},\vv{w}{})\in \PP$ and any $p,q\in Q$.

Any local rule $f$ induces a {\it global} rule $F$  as follows
\[
\left( \forall \bx\in X\right)\left(\forall \bi\in\C\right) \left(F(\bx)(\bi )=f(N_{\bx,\bi }) \right)\,.
\]
Having introduced the required notations, we now define a \emph{number-conserving} rule.
\begin{defn}
A local rule $f$ is called number-conserving if its corresponding global rule $F$ conserves the sum of states, \emph{i.e.}\ for each $\bx \in X$ it holds that
$\sigma(F(\bx))=\sigma(\bx)$.  
\end{defn}

\section{Number-conserving CAs for $d=2$}
In order to introduce our main ideas and rationale, we first describe the two-dimensional case. In Section 4, we will show how our arguments extend to arbitrary dimension.

\subsection{Basic properties}\label{sect:basic}

We first establish three simple, but useful facts. 

\begin{lem}\label{lemma:qq}
If a local rule $f$ is number-conserving, then each state is quiescent, \emph{i.e.} for any
$q\in Q$, it holds that
\begin{equation}\label{w1}
\ff{q}{q}{q}{q}{q} = q\,.
\end{equation}
\end{lem}

Proof. Consider a homogeneous configuration $\bq$ in which all cells have state $q$. Then $\sigma(\bq)=|\C|\cdot q$, while
% * <bernard.debaets@ugent.be> 2017-04-29T15:15:40.859Z:
%
% ^.
it is immediately clear that $\sigma\left( F(\bq)\right)= |\C|\cdot\ff{q}{q}{q}{q}{q}$, where $F$ is the corresponding global rule. Since $F$ conserves the sum of states, our statement follows.
\hfill\ensuremath{\square}

\begin{lem}\label{lemma:sumq}
If a local rule $f$ is number-conserving, then for any $q\in Q$, it holds that $\seff{q}{q}{q}{q}{q} = q$, \emph{i.e.}\
\begin{equation}\label{w2}
\ff{q}{0}{0}{0}{0}+\ff{0}{q}{0}{0}{0}+\ff{0}{0}{q}{0}{0}+\ff{0}{0}{0}{q}{0}+\ff{0}{0}{0}{0}{q} = q\,.
\end{equation}
\end{lem}
Proof. 
From Lemma \ref{lemma:qq}, we know that $\ff{0}{0}{0}{0}{0}=0$. Therefore, if we consider a configuration $\bx\in X$ in which only one cell has state $q$ and the other cells have state $0$, then $\sigma(\bx)= q$, while
it must hold at the same time that 
\[
\sigma\left( F(\bx) \right)= \ff{q}{0}{0}{0}{0}+\ff{0}{q}{0}{0}{0}+\ff{0}{0}{q}{0}{0}+\ff{0}{0}{0}{q}{0}+\ff{0}{0}{0}{0}{q}\,.
\]
So, \eqref{w2} is a direct consequence of the fact that
$F$ conserves the sum of states.
\hfill\ensuremath{\square}

The following lemma concerns the relationships between the values of the matching dimers. 

\begin{lem}\label{lemma:eq}
If a local rule $f$ is number-conserving, then for any
$q_1,q_2\in Q$, the following equalities hold
\begin{equation}\label{e1}
\ff{0}{q_1}{q_2}{0}{0}+\ff{0}{0}{q_1}{q_2}{0}=
\ff{0}{q_1}{0}{0}{0}+\ff{0}{0}{q_1}{0}{0}+\ff{0}{0}{q_2}{0}{0}+\ff{0}{0}{0}{q_2}{0}
\end{equation}
\begin{equation}\label{e2}
\ff{q_1}{0}{q_2}{0}{0}+\ff{0}{0}{q_1}{0}{q_2}=
\ff{q_1}{0}{0}{0}{0}+\ff{0}{0}{q_1}{0}{0}+\ff{0}{0}{q_2}{0}{0}+\ff{0}{0}{0}{0}{q_2}
\end{equation}
\begin{equation}\label{e3}
\ff{q_1}{0}{0}{q_2}{0}+\ff{0}{q_1}{0}{0}{q_2}=
\ff{q_1}{0}{0}{0}{0}+\ff{0}{q_1}{0}{0}{0}+\ff{0}{0}{0}{q_2}{0}+\ff{0}{0}{0}{0}{q_2}
\end{equation}
\begin{equation}\label{e4}
\ff{q_1}{q_2}{0}{0}{0}+\ff{0}{0}{0}{q_1}{q_2}=
\ff{q_1}{0}{0}{0}{0}+\ff{0}{0}{0}{q_1}{0}+\ff{0}{q_2}{0}{0}{0}+\ff{0}{0}{0}{0}{q_2}
\end{equation}
% * <bernard.debaets@ugent.be> 2017-04-29T15:30:47.454Z:
% 
% > \begin{equation}\label{e1}
% > \ff{0}{q_1}{q_2}{0}{0}+\ff{0}{0}{q_1}{q_2}{0}=
% > \ff{0}{q_1}{0}{0}{0}+\ff{0}{0}{q_1}{0}{0}+\ff{0}{0}{q_2}{0}{0}+\ff{0}{0}{0}{q_2}{0}
% > \end{equation}
% > \begin{equation}\label{e2}
% > \ff{q_1}{0}{q_2}{0}{0}+\ff{0}{0}{q_1}{0}{q_2}=
% > \ff{q_1}{0}{0}{0}{0}+\ff{0}{0}{q_1}{0}{0}+\ff{0}{0}{q_2}{0}{0}+\ff{0}{0}{0}{0}{q_2}
% > \end{equation}
% > \begin{equation}\label{e3}
% > \ff{q_1}{0}{0}{q_2}{0}+\ff{0}{q_1}{0}{0}{q_2}=
% > \ff{q_1}{0}{0}{0}{0}+\ff{0}{q_1}{0}{0}{0}+\ff{0}{0}{0}{q_2}{0}+\ff{0}{0}{0}{0}{q_2}
% > \end{equation}
% > \begin{equation}\label{e4}
% > \ff{q_1}{q_2}{0}{0}{0}+\ff{0}{0}{0}{q_1}{q_2}=
% > \ff{q_1}{0}{0}{0}{0}+\ff{0}{0}{0}{q_1}{0}+\ff{0}{q_2}{0}{0}{0}+\ff{0}{0}{0}{0}{q_2}
% > \end{equation}
% 
% ^.
\end{lem}
Proof. Each of the above equalities is easy to show if one considers a properly chosen configuration (see Fig.~\ref{configFIG}) and relies on Lemma~\ref{lemma:sumq}. Indeed, if we denote the configuration shown in 
Fig.~\ref{configFIG}(a) by $\bx$, then from  $\sigma(\bx)=\sigma\left( F(\bx) \right)$, it follows that
\begin{equation}\label{fa1}
\begin{aligned}
q_1+q_2 = {} &\ff{0}{0}{0}{0}{q_1}+\ff{0}{0}{0}{0}{q_2}+\ff{0}{0}{0}{q_1}{0}+\ff{0}{0}{q_1}{q_2}{0}\\
&+\ff{0}{q_1}{q_2}{0}{0}+\ff{0}{q_2}{0}{0}{0}+\ff{q_1}{0}{0}{0}{0}+\ff{q_2}{0}{0}{0}{0}.
\end{aligned}
\end{equation}
From Lemma \ref{lemma:sumq}, we know that
\begin{equation}\label{fa2}
q_1-\ff{0}{0}{0}{0}{q_1}-\ff{0}{0}{0}{q_1}{0}-\ff{q_1}{0}{0}{0}{0}=\ff{0}{q_1}{0}{0}{0}+\ff{0}{0}{q_1}{0}{0}
\end{equation}
and
\begin{equation}\label{fa3}
q_2-\ff{0}{0}{0}{0}{q_2}-\ff{0}{q_2}{0}{0}{0}-\ff{q_2}{0}{0}{0}{0}=\ff{0}{0}{q_2}{0}{0}+\ff{0}{0}{0}{q_2}{0}.
\end{equation}
Combining Eqs.\ (\ref{fa1})--(\ref{fa3}), our claim follows easily.
The proofs of Eqs.~(\ref{e2})--(\ref{e4}) are similar, and are based on the configurations in Figs.~\ref{configFIG}(b)--(d).
\hfill\ensuremath{\square}

One can see that the terms in the right-hand sides of Eqs.\ (\ref{e1})--(\ref{e4}) are the monomer expansions of the terms in the left-hand sides. 

\begin{figure}[!h]
\begin{center}
 \subfloat[]{
 \begin{tabular}{|c|c|c|c|} 
 \hline
 0 & 0 & 0 & 0 \\
 \hline
 0 & $q_1$ & $q_2$ & 0\\
 \hline
 0 & 0 & 0 & 0 \\
 \hline
 0 & 0 & 0 & 0\\ 
\hline
\end{tabular}
 }
 \quad
 \subfloat[]{
 \begin{tabular}{|c|c|c|c|}
 \hline
 0 & 0 & 0 \\
 \hline
 0 & $q_1$ & 0 \\
 \hline
 0 & $q_2$ & 0 \\
 \hline
 0 & 0 & 0\\
 \hline
\end{tabular}
 }
 \quad
 \subfloat[]{
 \begin{tabular}{|c|c|c|c|}
 \hline
 0 & 0 & 0 & 0 \\
 \hline
 0 & $q_1$ & 0 & 0\\
 \hline
 0 & 0 & $q_2$ & 0 \\
 \hline
 0 & 0 & 0 & 0\\ 
 \hline
 \end{tabular}
 }
 \quad
 \subfloat[]{
 \begin{tabular}{|c|c|c|c|} 
 \hline
 0 & 0 & 0 & 0 \\
 \hline
 0 & 0 & $q_1$ & 0\\
 \hline
 0 & $q_2$ & 0 & 0 \\
 \hline
 0 & 0 & 0 & 0 \\
 \hline
\end{tabular}
 }
\end{center} 
\caption{The configurations used to show the relationships between the matching dimers in Lemma~ \ref{lemma:eq}: (a)~horizontal, (b)~vertical, (c)--(d) diagonal arrangements.}
\label{configFIG}
 \end{figure}

\subsection{Necessary and sufficient conditions}

Using the lemmas presented in Section \ref{sect:basic}, we can prove the main theorem for $d=2$.
\begin{thm}\label{main}
A local rule $f$ is number-conserving if and only if 
for any $q_1, q_2, q_3, q_4, q_5\in Q$, it holds that
\begin{equation}
    \label{dc}
\begin{aligned}
    \ff{q_1}{q_2}{q_3}{q_4}{q_5} = & q_1 +
    \ff{0}{q_2}{0}{0}{q_5} - \ff{0}{q_1}{0}{0}{q_4}+
    \ff{0}{0}{0}{q_4}{q_5} - \ff{0}{0}{0}{q_1}{q_2} \\
     &+ \ff{0}{0}{q_3}{q_4}{0} - \ff{0}{0}{q_2}{q_3}{0} 
     + \ff{0}{0}{q_3}{0}{q_5} - \ff{0}{0}{q_1}{0}{q_3}  \\
     &+\ff{0}{0}{q_2}{0}{0} - \ff{0}{0}{q_3}{0}{0}
     +\ff{0}{0}{0}{q_3}{0} - \ff{0}{0}{0}{q_4}{0} \\
     &+ \ff{0}{0}{0}{0}{q_2} +\ff{0}{0}{0}{0}{q_3} 
     +\ff{0}{0}{0}{0}{q_4} -
     \ff{0}{0}{0}{0}{q_1} -
     2\ff{0}{0}{0}{0}{q_5}\,.
\end{aligned}
\end{equation}
\,
\end{thm}

Proof: To show that \eqref{dc} is necessary, let us consider the configuration in Fig.~\ref{fig:q1q5}. Suppose that $f$ is number-conserving, then it holds that
\begin{equation}
    \label{r1}
\begin{aligned}
    q_1+q_2+q_3+q_4+q_5=\ff{0}{0}{0}{0}{q_1}+\ff{0}{0}{0}{q_1}{q_2}+\ff{0}{0}{q_1}{0}{q_3}\\
    +\ff{0}{q_1}{0}{0}{q_4}+\ff{0}{0}{0}{q_2}{0}+\ff{0}{0}{q_2}{q_3}{0}+\ff{q_1}{q_2}{q_3}{q_4}{q_5}+\ff{0}{q_3}{q_4}{0}{0}\\
    +\ff{0}{q_4}{0}{0}{0}+\ff{q_2}{0}{0}{q_5}{0}+\ff{q_3}{0}{q_5}{0}{0}+\ff{q_4}{q_5}{0}{0}{0}+
    \ff{q_5}{0}{0}{0}{0}\,.
\end{aligned}
\end{equation}
As expected, in the right-hand side of \eqref{r1} we retrieve eight dimers, because there are exactly eight pairs in $\PP$. Now, for every two matching pairs we select one such that
\[\CP=\{ 
\lbrace \vv{0}{},\vv{v}{1} \rbrace, 
\lbrace \vv{0}{},\mbox{-}\vv{v}{2}\rbrace, 
\lbrace \mbox{-}\vv{v}{1},\mbox{-}\vv{v}{2}\rbrace, \lbrace \vv{v}{1},\mbox{-}\vv{v}{2}\rbrace \}\,.
\]
Using Lemma~\ref{lemma:eq}, we represent the dimers for the pairs not present in $\CP$ in terms of their matching dimers:
\begin{subequations}
% * <bernard.debaets@ugent.be> 2017-04-29T15:44:54.983Z:
% 
% use eqnnarray
% \arraycolsep=2pt
% \begin{eqnarray}
%    &=&
% ...
% to align better on =, and equal spaces around =
% 
% 
% ^.
\begin{align}
\ff{0}{q_3}{q_4}{0}{0}=&
\ff{0}{q_3}{0}{0}{0}+\ff{0}{0}{q_3}{0}{0}+\ff{0}{0}{q_4}{0}{0}+\ff{0}{0}{0}{q_4}{0}-\ff{0}{0}{q_3}{q_4}{0}      \label{r2:a}
\\
\ff{q_3}{0}{q_5}{0}{0}=&
\ff{q_3}{0}{0}{0}{0}+\ff{0}{0}{q_3}{0}{0}+\ff{0}{0}{q_5}{0}{0}+\ff{0}{0}{0}{0}{q_5}-\ff{0}{0}{q_3}{0}{q_5} \label{r2:b}\\
\ff{q_2}{0}{0}{q_5}{0}=&
\ff{q_2}{0}{0}{0}{0}+\ff{0}{q_2}{0}{0}{0}+\ff{0}{0}{0}{q_5}{0}+\ff{0}{0}{0}{0}{q_5}-\ff{0}{q_2}{0}{0}{q_5} \label{r2:c}\\
\ff{q_4}{q_5}{0}{0}{0}=&
\ff{q_4}{0}{0}{0}{0}+\ff{0}{0}{0}{q_4}{0}+\ff{0}{q_5}{0}{0}{0}+\ff{0}{0}{0}{0}{q_5}-\ff{0}{0}{0}{q_4}{q_5}\,.\label{r2:d}
\end{align}
\end{subequations}
Combining \eqref{r1} and Eqs.\ (\ref{r2:a})--(\ref{r2:d}) and applying Lemma \ref{lemma:sumq} for $q=q_i$, $i\in\{ 2,3,4,5\}$, we obtain \eqref{dc}.

Conversely, one can easily see that \eqref{dc} is sufficient, but we give a detailed proof for the sake of completeness. Let us assume that \eqref{dc} holds and let $\bx$ be any configuration. Then it holds that
\begin{equation}
\begin{aligned}
& \sum_{(i,j)\in \C}\ff{x_{i,j-1}}{x_{i-1,j}}{x_{i,j}}{x_{i+1,j}}{x_{i,j+1}} =\sum_{(i,j)\in \C}x_{i,j-1} +
\sum_{(i,j)\in \C}\left[\ff{0}{x_{i-1,j}}{0}{0}{x_{i,j+1}} -\ff{0}{x_{i,j-1}}{0}{0}{x_{i+1,j}}\right]  \\
+&\sum_{(i,j)\in \C} \left[\ff{0}{0}{0}{x_{i+1,j}}{x_{i,j+1}} - \ff{0}{0}{0}{x_{i,j-1}}{x_{i-1,j}} \right] 
+ \sum_{(i,j)\in \C} \left[\ff{0}{0}{x_{i,j}}{x_{i+1,j}}{0} - \ff{0}{0}{x_{i-1,j}}{x_{i,j}}{0}\right]    \\
+& 
\sum_{(i,j)\in \C} \left[\ff{0}{0}{x_{i,j}}{0}{x_{i,j+1}} - \ff{0}{0}{x_{i,j-1}}{0}{x_{i,j}}\right] 
+\sum_{(i,j)\in \C}\left[\ff{0}{0}{x_{i-1,j}}{0}{0} -\ff{0}{0}{x_{i,j}}{0}{0}\right]\\
+& \sum_{(i,j)\in \C}\left[\ff{0}{0}{0}{x_{i,j}}{0}- \ff{0}{0}{0}{x_{i+1,j}}{0} \right] 
+ \sum_{(i,j)\in \C}\ff{0}{0}{0}{0}{x_{i-1,j}} 
+ \sum_{(i,j)\in \C}\ff{0}{0}{0}{0}{x_{i,j}} \\
+& \sum_{(i,j)\in \C}\ff{0}{0}{0}{0}{x_{i+1,j}} - \sum_{(i,j)\in \C}\ff{0}{0}{0}{0}{x_{i,j-1}}  - 2\sum_{(i,j)\in \C}\ff{0}{0}{0}{0}{x_{i,j+1}}\,.
\end{aligned}
\label{eqX}
\end{equation}
It is obvious that
\[
\sum_{(i,j)\in \C}\ff{0}{x_{i-1,j}}{0}{0}{x_{i,j+1}}=\sum_{(i-1,j+1)\in \C}\ff{0}{x_{i,j-1}}{0}{0}{x_{i+1,j}}=\sum_{(i,j)\in \C}\ff{0}{x_{i,j-1}}{0}{0}{x_{i+1,j}}\,,
\]
assuming periodic boundary conditions. Thus, we get that
\[
\sum_{(i,j)\in \C}\left[\ff{0}{x_{i-1,j}}{0}{0}{x_{i,j+1}} -\ff{0}{x_{i,j-1}}{0}{0}{x_{i+1,j}}\right]=0\,,
\]
and likewise all terms of type $ \sum\left[\ldots -\ldots \right]$ are equal to $0$. Moreover, following the same reasoning, we see that
\begin{equation*}
\begin{aligned}
\sum_{(i,j)\in \C}\ff{0}{0}{0}{0}{x_{i-1,j}}=\sum_{(i,j)\in \C}\ff{0}{0}{0}{0}{x_{i,j}}=
\sum_{(i,j)\in C}\ff{0}{0}{0}{0}{x_{i+1,j}}\\
= \sum_{(i,j)\in \C}\ff{0}{0}{0}{0}{x_{i,j-1}} =\sum_{(i,j)\in \C}\ff{0}{0}{0}{0}{x_{i,j+1}}\,.
\end{aligned}
\end{equation*}
Using these equalities in \eqref{eqX}, we obtain
\begin{equation}
\label{eqY}
\sum_{(i,j)\in \C}\ff{x_{i,j-1}}{x_{i-1,j}}{x_{i,j}}{x_{i+1,j}}{x_{i,j+1}} = \sum_{(i,j)\in \C}x_{i,j-1}=
\sum_{(i,j)\in \C}x_{i,j}\,,
\end{equation}
which implies that $f$ is number-conserving.
\hfill\ensuremath{\square}

\begin{figure}
\begin{center}
 \begin{tabular}{|c|c|c|c|c|}
 \hline
 0 & 0 & 0 & 0 & 0\\
 \hline
 0 & 0 & $q_1$ & 0 & 0\\
 \hline
 0 & $q_2$ & $q_3$ & $q_4$ & 0 \\
 \hline
 0 & 0 & $q_5$ & 0 & 0 \\ 
 \hline
 0 & 0 & 0 & 0 & 0\\
 \hline
\end{tabular}
\end{center}
\caption{The essential part of the configuration used in the proof of Theorem \ref{main}. All remaining cells of the grid $\C$ have state $0$.}
\label{fig:q1q5}
\end{figure}

We remark that one can formulate Theorem \ref{main} in at least different $5\cdot 2^4$ ways. Firstly, one should choose an entry $q_i$ as the leading term, after which one has to select terms with two nonzero entries to expand using Lemma \ref{lemma:eq}. Then, one has to fill the right-hand side with monomers to match the monomer expansion of the left-hand side.
An example of an alternative formula with $q_3$ as the leading term reads:

\begin{equation}
    \label{dc2}
\begin{aligned}
    \ff{q_1}{q_2}{q_3}{q_4}{q_5} =& q_3 +
    \ff{q_1}{q_2}{0}{0}{0} - \ff{q_4}{q_5}{0}{0}{0}+
    \ff{0}{q_2}{0}{0}{q_5} - \ff{0}{q_1}{0}{0}{q_4} \\
     &+ \ff{0}{0}{q_3}{q_4}{0} - \ff{0}{0}{q_2}{q_3}{0} 
     +  \ff{q_1}{0}{q_3}{0}{0} - \ff{q_3}{0}{q_5}{0}{0}  \\
     &+\ff{q_4}{0}{0}{0}{0} - \ff{q_1}{0}{0}{0}{0}
     +\ff{0}{0}{0}{0}{q_4} - \ff{0}{0}{0}{0}{q_3} \\
     &+ \ff{0}{q_1}{0}{0}{0} + \ff{0}{q_5}{0}{0}{0}- \ff{0}{q_2}{0}{0}{0}- \ff{0}{q_3}{0}{0}{0} \\
% * <bernard.debaets@ugent.be> 2017-04-29T15:47:58.099Z:
%
% ^.
     & +\ff{0}{0}{q_2}{0}{0} +\ff{0}{0}{q_5}{0}{0} -  2\ff{0}{0}{q_3}{0}{0}\,.
\end{aligned}
\end{equation}
 One can see that this expression contains more terms than~\eqref{dc}.

\begin{cor}
To define a number-conserving rule, it is sufficient to specify its values for nontrivial monomers and dimers for pairs belonging to $\CP$. If for any $q_1, q_2, q_3, q_4, q_5\in Q$, the value of the right-hand side of \eqref{dc} belongs to $Q$, then this rule is number-conserving.
\label{cor:check}
\end{cor}

Using Corollary \ref{cor:check} it is very easy to enumerate all  number-conserving two-dimensional CAs in the case of few states. For example, if we consider $Q=\{ 0,1\}$, then there are only nine different number-conserving rules, namely the identity rule and the shift and traffic rules in each of the four possible directions. This means that they are the two-dimensional versions of one-dimensional rules in the horizontal or vertical direction. 
The case $Q=\{0,1,2\}$ is much more interesting.
There are $1327$ number-conserving rules and only 287 of them are extensions of horizontal and vertical one-dimensional ones. 

\section{The case of an arbitrary dimension $d$}

\subsection{Basic properties}

In this section, we consider an arbitrary dimension $d$ and we establish necessary and sufficient conditions for CAs with the von Neumann neighborhood to be number-conserving. 
We formulate the conditions in terms of the local rule of the given CA.

Analogously to the case $d=2$, we start by showing the following lemma, which is the $d$-dimensional counterpart of Lemmata \ref{lemma:qq}  and \ref{lemma:sumq}.
\begin{lem}\label{lemma:qq3}
If a local rule $f$ is number-conserving, then each state is quiescent, i.e. for any
$q\in Q$, it holds that $f(H_q)=q$. Moreover, $f^{E}(H_q)=q$.
\end{lem}

Proof: The first statement is well known and can be shown by considering a homogeneous configuration $\bx$, in which every cell has state $q$. To prove the second one, let us consider a~configuration~$\bx$ in which only one cell has state $q$ and the other ones have state $0$:
\[
\bx (\bi )=
\begin{cases}
q  \quad\text{, if }  \bi = \bze\\
0  \quad\text{, if }  \bi \neq \bze.
\end{cases} 
\]
Obviously, $\sigma (\bx )=q$. Note that if $\bze\not\in P(\bi )$, then the neighborhood configuration $N_{\bx,\bi}$ is trivial, while if $\bze\in P(\bi )$, then the neighborhood configuration $N_{\bx,\bi}$ is a monomer. Moreover, if $P(\bi )$ contains~$\bze$, then it holds that $\bze =\bi+\vv{v}{}$ for some $\vv{v}{}\in V$ and then $N_{\bx,\bi}(\vv{v}{})=N_{\bx,\bze-\vv{v}{} }(\vv{v}{})=\bx (\bze )=q$, \emph{i.e.}\ $N_{\bx,\bze-\vv{v}{} }=M_{\vv{v}{}:q}$. Thus, we have
\[
\sigma (F(\bx))=\sum_{\bi\in\C}f(N_{\bx,\bi})=\sum_{\vv{v}{}\in V}f(N_{\bx,\bze-\vv{v}{} })=
\sum_{\vv{v}{}\in V}f(M_{\vv{v}{}:q})=f^{E}(H_q).
\]
As  $\sigma (\bx )=\sigma (F(\bx))$, we obtain our claim.

The following fact concerns properties of pairs of matching dimers. For $d=2$, we listed all four equations in Lemma~\ref{lemma:eq}. Here we write $d^2$ equations as a single one using the dimer notation.

\begin{lem}\label{lemma:cp3}
If a local rule $f$ is number-conserving, then for any $(\vv{u}{},\vv{w}{})\in \PP$ and any $p,q\in Q$, it holds that
\begin{equation}\label{DD}
f\left( 
\dimer{\vv{u}{}:p}{
\vv{w}{}:q}
\right) 
+f\left(\dimer{\mbox{-}\vv{w}{}:p}{
\mbox{-}\vv{u}{}:q}\right) =
f^{E}\!\left( \dimer{\vv{u}{}:p}{
\vv{w}{}:q}\right) +f^{E}\left( \dimer{\mbox{-}\vv{w}{}:p}{
\mbox{-}\vv{u}{}:q}\right)\,.
\end{equation}
\end{lem}

Proof: Let $\lbrace \vv{u}{},\vv{w}{}\rbrace\in \PP$ and $p,q\in Q$. Further, we consider the following configuration $\bx$:
\[
\bx (\bi )=
\begin{cases}
p \quad \text{, if } \bi = \bze+\vv{u}{}\\
q \quad\text{, if } \bi = \bze+\vv{w}{}\\
0 \quad \text{, otherwise.}\\
\end{cases}
\]
Clearly, $\sigma (\bx )=p+q$. If neither $\bze+\vv{u}{}$, nor $\bze+\vv{w}{}$ belongs to $P(\bi)$, then $N_{\bx,\bi }$ is trivial, so it suffices to consider the following 3 cases:
\begin{itemize}
\item[(a)] $\bze+\vv{u}{}\in P(\bi )$ and $\bze+\vv{w}{}\not\in P(\bi )$: in this case, $N_{\bx,\bi }$ is a monomer,
\item[(b)] $\bze+\vv{u}{}\not\in P(\bi )$ and $\bze+\vv{w}{}\in P(\bi )$: also in this case, $N_{\bx,\bi }$ is a monomer,
\item[(c)] $\bze+\vv{u}{}\in P(\bi )$ and $\bze+\vv{w}{}\in P(\bi )$: in this case, $N_{\bx,\bi }$ is a dimer.
\end{itemize}
Let us recall that $P(\bi)$ contains $\bze+\vv{u}{}$ if and only if $\bze+\vv{u}{} =\bi +\vv{v}{}$ for some $\vv{v}{}\in V$, and more importantly, it then holds that $N_{\bx,\bi}(\vv{v}{})=N_{\bx,\bze + \vv{u}{} - \vv{v}{} }(\vv{v}{})=\bx (\bze+\vv{u}{})$. 
Similarly, $P(\bi)$ contains $\bze+\vv{w}{}$ if and only if $\bze+\vv{w}{} =\bi +\vv{v}{}'$ for some $\vv{v}{}'\in V$ and then $N_{\bx,\bi } (\vv{v}{}')=N_{\bx,\bze+\vv{w}{}-\vv{v}{}'}(\vv{v}{}')=\bx (\bze+\vv{w}{})$. 
Thus $P(\bi )$ contains both $\bze+\vv{u}{}$ and $\bze+\vv{w}{}$ only when  $\bi = \bze + \vv{u}{} - \vv{v}{}=
 \bze + \vv{w}{} - \vv{v}{}'$ for some $\vv{v}{},\vv{v}{}'\in V$. As $\lbrace \vv{u}{},\vv{w}{}\rbrace \in \PP$ then $\vv{u}{}$ and $\vv{w}{}$ 
act on different components and at least one of them is nonzero. Hence, the vector equation 
$\bi = \bze + \vv{u}{} - \vv{v}{}= \bze + \vv{w}{} - \vv{v}{}'$
 has exactly two solutions:
\[
(\vv{v}{},\vv{v}{}')=(\mbox{-}\vv{w}{},\mbox{-}\vv{u}{}) \quad\mbox{and}\quad
(\vv{v}{},\vv{v}{}')=(\vv{u}{},\vv{w}{})\,. 
\]
From this, we conclude that $P(\bi )$ contains both $\bze+\vv{u}{}$ and $\bze+\vv{w}{}$ only when $\bi=\bze$ or 
$\bi=\bze+\vv{u}{}+\vv{w}{}$. 
Moreover, 
\begin{equation}\label{dimers}
N_{\bx,\bze }=\dimer{\vv{u}{}:\bx (\bze+\vv{u}{})}{\vv{w}{}:\bx (\bze+\vv{w}{})} =\dimer{\vv{u}{}:p}{
\vv{w}{}:q} \quad\text{ and }\quad N_{\bx,\bze+\vv{u}{}+\vv{w}{} }=
\dimer{\mbox{-}\vv{w}{}:\bx (\bze+\vv{u}{})}{\mbox{-}\vv{u}{}:\bx (\bze+\vv{w}{})}=\dimer{\mbox{-}\vv{w}{}:p}{
\mbox{-}\vv{u}{}:q}\,.
\end{equation}
Summarizing the above observations, we get
\begin{itemize}
\item[(a)] $\bze+\vv{u}{}\in P(\bi )$ and $\bze+\vv{w}{}\not\in P(\bi )$ if and only if $\bi = \bze+\vv{u}{}-\vv{v}{}$ where $\vv{v}{}\in V\setminus \{\vv{u}{},\mbox{-}\vv{w}{}\}$ and in this case $N_{\bx,\bi }= M_{\vv{v}{}:p}$;
\item[(b)] $\bze+\vv{u}{}\not\in P(\bi )$ and $\bze+\vv{w}{}\in P(\bi )$ if and only if $\bi = \bze+\vv{w}{}-\vv{v}{}'$ where $\vv{v}{}'\in V\setminus \{\vv{w}{},\mbox{-}\vv{u}{}\}$ and in this case $N_{\bx,\bi }= M_{\vv{v}{}':q}$;
\item[(c)] $\bze+\vv{u}{}\in P(\bi )$ and $\bze+\vv{w}{}\in P(\bi )$ if and only if $\bi = \bze $ or 
$\bi=\bze+\vv{u}{}+\vv{w}{}$ and in this case $N_{\bx,\bi }$ is given by~\eqref{dimers}.
\end{itemize}

Now, we can write
\begin{equation}
\label{stareq}
\begin{aligned}
\sigma(F(\bx )) &= \sum_{\bi\in\C}f(N_{\bx,\bi })
\\
&=\sum_{\vv{v}{}\in V\setminus\{ \vv{u}{},\mbox{-}\vv{w}{}\}}\hspace{-6mm}f(N_{\bx,\bze+\vv{u}{}-\vv{v}{} }) 
+ \sum_{\vv{v}{}'\in V\setminus\{ \vv{w}{},\mbox{-}\vv{u}{}\}}\hspace{-6mm}f(N_{\bx,\bze+\vv{w}{}-\vv{v}{}'})+
 f(N_{\bx,\bze })+f(N_{\bx,\bze+\vv{u}{}+\vv{w}{} })
\\
&=\sum_{\vv{v}{}\in V\setminus\{ \vv{u}{},\mbox{-}\vv{w}{}\}}\hspace{-6mm}f\left( M_{\vv{v}{}:p}\right)+\sum_{\vv{v}{}'\in V\setminus\{\vv{w}{}, \mbox{-}\vv{u}{}\}}\hspace{-6mm} f\left( M_{\vv{v}{}':q}\right)+f\left( 
\dimer{\vv{u}{}:p}{
\vv{w}{}:q}
\right) 
+f\left(\dimer{\mbox{-}\vv{w}{}:p}{
\mbox{-}\vv{u}{}:q}\right)\,.
\end{aligned}
\end{equation}
From Lemma~\ref{lemma:qq3}, we know that 
$\displaystyle \sum_{\vv{v}{}\in V}f\left( M_{\vv{v}{}:p}\right) = 
f^{E}(H_p)=p$, thus 
\begin{equation}
\label{eq21}
\begin{aligned}
\sum_{\vv{v}{}\in V\setminus\{ \vv{u}{},\mbox{-}\vv{w}{}\}}
\!\!\!\!\! f\left( M_{\vv{v}{}:p}\right) &= \sum_{\vv{v}{}\in V}f\left( M_{\vv{v}{}:p}\right)-
\left[ f\left( M_{\vv{u}{}:p}\right)+f\left( M_{\mbox{-}\vv{w}{}:p}\right)\right]\\
&= p-f\left( M_{\vv{u}{}:p}\right)-f\left( M_{\mbox{-}\vv{w}{}:p}\right)
\end{aligned}
\end{equation}
and, similarly,
\begin{equation}
\label{eq22}
\sum_{\vv{v}{}'\in V\setminus\{ \vv{w}{},\mbox{-}\vv{u}{}\}} f\left( M_{\vv{v}{}':q}\right)=
q-f\left( M_{\vv{w}{}:q}\right)-f\left( M_{\mbox{-}\vv{u}{}:q}\right)\,.
\end{equation}
Combining Eqs. (\ref{stareq})--(\ref{eq22}), we obtain
\begin{equation}
\begin{aligned}
\sigma(F(\bx ))= &p+q- f\left( M_{\vv{u}{}:p}\right)-f\left( M_{\mbox{-}\vv{w}{}:p}\right)-f\left( M_{\vv{w}{}:q}\right)-f\left( M_{\mbox{-}\vv{u}{}:q}\right)+f\left( 
\dimer{\vv{u}{}:p}{
\vv{w}{}:q}
\right) 
+f\left(\dimer{\mbox{-}\vv{w}{}:p}{
\mbox{-}\vv{u}{}:q}\right) \\
=  p+q- & \left( f\left( M_{\vv{u}{}:p}\right)+f\left( M_{\vv{w}{}:q}\right) \right)-\left( f\left( M_{\mbox{-}\vv{w}{}:p}\right)+ f\left( M_{\mbox{-}\vv{u}{}:q}\right)\right) +f\left( 
\dimer{\vv{u}{}:p}{
\vv{w}{}:q}
\right) 
+f\left(\dimer{\mbox{-}\vv{w}{}:p}{
\mbox{-}\vv{u}{}:q}\right) \\
= p+q- & f^{E}\left( \dimer{\vv{u}{}:p}{
\vv{w}{}:q}\right)-f^{E}\left( \dimer{\mbox{-}\vv{w}{}:p}{
\mbox{-}\vv{u}{}:q}\right) +f\left( 
\dimer{\vv{u}{}:p}{
\vv{w}{}:q}
\right) 
+f\left(\dimer{\mbox{-}\vv{w}{}:p}{
\mbox{-}\vv{u}{}:q}\right)\,,
\end{aligned}
\end{equation}
where we used the fact that the monomer expansion of a dimer is given by \eqref{dimerexp}. As $\sigma (\bx )=\sigma (F(\bx ))=p+q$, the above 
collapses to \eqref{DD}.
\hfill\ensuremath{\square}

\subsection{Necessary and sufficient condition}
The following theorem presents one of $(2d+1)\cdot 2^{d^2}$ possible formulations of necessary and sufficient conditions for a local rule to be number-conserving. This is because we have $2d+1$ possibilities for choosing the leading term ($\vv{\eta}{}$ below) and $2^{d^2}$ possibilities for choosing $\CP$. 
\begin{thm}\label{maind}
Let $\vv{\eta}{}\in V$ and $\CP$ be fixed. A local rule $f$ is number-conserving if and only if for any
$N\in\NN$, it holds that
\begin{equation}
\label{eq:thm43}
\begin{aligned}
f(N) =& N(\vv{\eta}{})+\sum_{(\vv{u}{},\vv{w}{})\in \CP}\left[ 
f\left( \dimer{\vv{u}{}:N(\vv{u}{})}
{\vv{w}{}:N(\vv{w}{})}\right)
 -
f\left(  \dimer{\vv{u}{}:N(\mbox{-}\vv{w}{})}
{\vv{w}{}:N(\mbox{-}\vv{u}{})}\right)\right] 
+\sum_{\vv{v}{}\in V\setminus\{\vv{\eta}{}\}}f^{E}\left( H_{N(\vv{v}{})}\right)
\\
&-\sum_{(\vv{u}{},\vv{w}{})\in \CP}
\left[f^{E}\left( \dimer{\vv{u}{}:N(\vv{u}{})}
{\vv{w}{}:N(\vv{w}{})}\right)+
f^{E}\left( \dimer{\mbox{-}\vv{w}{}:N(\vv{u}{})}
{\mbox{-}\vv{u}{}:N(\vv{w}{})}\right)\right]
-\sum_{\vv{v}{}\in V_+}
f\left( M_{\vv{v}{}:N(\mbox{-}\vv{v}{})} \right)\,.
\end{aligned}
\end{equation}
\end{thm}
Proof: 
It is easy to see that this condition is sufficient. Indeed, consider the grid $\C$ and an arbitrary configuration $\bx$. Now, it is enough to fix $\vv{\eta}{}\in V$, write \eqref{eq:thm43} for each neighborhood configuration $N_{\bx,\bi}$ and combine the expressions. Doing so, on the left-hand side, we obtain  $\sigma (F(\bx ))$, while on the right-hand side the dimers cancel out analogously to the case $d=2$ (see how~\eqref{eqX} implies~\eqref{eqY}). 
For monomers, it suffices to check that their number with a positive sign in
$\displaystyle \sum_{\vv{v}{}\in V\setminus\{\vv{\eta}{}\}}f^{E}\left( H_{N(\vv{v}{})}\right)$, 
which equals $2d(2d+1)$, equals the number of monomers with a minus sign, which equals $d^2\cdot(2+2) +2d$.
Hence, we obtain $\sigma (\bx )$ on the right-hand side.  

To prove that \eqref{eq:thm43} is necessary, let us fix $N\in \NN$ and consider the configuration $\bx$ defined by
\begin{equation*}
\bx (\bi )=
\begin{cases}
N(\vv{v}{})  &\quad\text{, if }\ \bi = \bze +\vv{v}{} \text{ for some } \vv{v}{}\in V,\\
0            &\quad\text{, otherwise.}
\end{cases}
\end{equation*}
Note that $N_{\bx ,\bze }=N$. By Lemma \ref{lemma:qq3}, it follows that for any $\vv{v}{}\in V$, it holds that $N(\vv{v}{})=f^{E}\left( H_{N(\vv{v}{})}\right)$, so
\[
\sigma (\bx )=\sum_{\vv{v}{}\in V}N(\vv{v}{})= N(\vv{\eta}{})+\sum_{\vv{v}{}\in V\setminus\{\vv{\eta}{}\}} N(\vv{v}{})
= N(\vv{\eta}{})+\sum_{\vv{v}{}\in V\setminus\{\vv{\eta}{}\}}f^{E}\left( H_{N(\vv{v}{})}\right)\,.
\]
On the other hand, as configuration $\bx$ is zero outside $P(\bze)$, from Lemma \ref{lemma:vN} we obtain
\begin{equation}
\label{eqSigmaF}
\sigma (F(\bx ))=\sum_{\bi\in\C} f(N_{\bx,\bi })=f(N_{\bx,\bze})+
\sum_{\{\bi \mid\; |P(\bi )\cap P(\bze )|=2\}}\hspace{-6mm} f(N_{\bx,\bi })+\sum_{\{\bi \mid\; |P(\bi )\cap P(\bze )|=1\}}\hspace{-6mm} f(N_{\bx,\bi })\,.
\end{equation}
In view of Remark \ref{remark}, if $|P(\bi )\cap P(\bze )|=2$, then there exists a unique pair $(\vv{u}{},\vv{w}{})\in \PP$ such that $\bi =\bze + \vv{u}{} + \vv{w}{}$ and $P(\bi )\cap P(\bze )=\{\bze + \vv{u}{},\bze + \vv{w}{}\}$. Hence, $N_{\bx,\bi }$ is a dimer satisfying
\begin{equation}
\label{dimerxi}
N_{\bx,\bze + \vv{u}{} + \vv{w}{}}(\mbox{-}\vv{u}{})=\bx (\bze + \vv{w}{})=N(\vv{w}{})
\quad \text{ and }\quad
N_{\bx,\bze + \vv{u}{} + \vv{w}{}}(\mbox{-}\vv{w}{})=\bx (\bze + \vv{u}{})=N(\vv{u}{}),
\end{equation}
which means that $N_{\bx,\bze + \vv{u}{} + \vv{w}{}}=
\dimer{\mbox{-}\vv{w}{}:N(\vv{u}{})}{\mbox{-}\vv{u}{}:N(\vv{w}{})}$. Hence, for a given $\CP$, we have
\begin{equation}
\label{eqFromRk22}
\sum_{\{\bi \mid\; |P(\bi )\cap P(\bze )|=2\}} \hspace{-6mm} f(N_{\bx,\bi })=\hspace{-3mm} \sum_{(\vv{u}{},\vv{w}{})\in \PP} \hspace{-1mm} f\left( \dimer{\mbox{-}\vv{w}{}:N(\vv{u}{})}{\mbox{-}\vv{u}{}:N(\vv{w}{})}
\right)=
\sum_{(\vv{u}{},\vv{w}{})\in \CP}\left[  f\left( \dimer
{\mbox{-}\vv{w}{}:N(\vv{u}{})}{\mbox{-}\vv{u}{}:N(\vv{w}{})}\right) + f\left( \dimer
{\vv{u}{}:N(\mbox{-}\vv{w}{})}{\vv{w}{}:N(\mbox{-}\vv{u}{})}\right)\right]\,,
\end{equation}
while from Lemma~\ref{lemma:cp3}, we have that for any $(\vv{u}{},\vv{w}{})\in \PP$:
\begin{equation}
\label{eqFromRk42}
f\left( \dimer{\mbox{-}\vv{w}{}:N(\vv{u}{})}{\mbox{-}\vv{u}{}:N(\vv{w}{})}\right)=
f^{E}\left( \dimer{\vv{u}{}:N(\vv{u}{})}
{\vv{w}{}:N(\vv{w}{})}\right)+
f^{E}\left( \dimer{\mbox{-}\vv{w}{}:N(\vv{u}{})}
{\mbox{-}\vv{u}{}:N(\vv{w}{})}\right)-
f\left( \dimer{\vv{u}{}:N(\vv{u}{})}
{\vv{w}{}:N(\vv{w}{})}\right) .
\end{equation}
and from Eqs.~(\ref{eqFromRk22}) and (\ref{eqFromRk42}), we conclude that
\begin{equation}
\label{eqDot}
\begin{aligned}
\sum_{\{\bi \mid\; |P(\bi )\cap P(\bze )|=2\}} \hspace{-3mm} f(N(\bx,\bi )) =&
\sum_{(\vv{u}{},\vv{w}{})\in \CP}\left[  f\left( \dimer
{\vv{u}{}:N(\mbox{-}\vv{w}{})}{\vv{w}{}:N(\mbox{-}\vv{u}{})}\right) - 
f\left( \dimer{\vv{u}{}:N(\vv{u}{})}
{\vv{w}{}:N(\vv{w}{})}\right)\right] \\
&+ \sum_{(\vv{u}{},\vv{w}{})\in \CP}
\left[f^{E}\left( \dimer{\vv{u}{}:N(\vv{u}{})}
{\vv{w}{}:N(\vv{w}{})}\right)+
f^{E}\left( \dimer{\mbox{-}\vv{w}{}:N(\vv{u}{})}
{\mbox{-}\vv{u}{}:N(\vv{w}{})}\right)\right].
\end{aligned}
\end{equation}

Again using Lemma \ref{lemma:vN}, we know that $|P(\bi )\cap P(\bze )|=1$ only when $\bi =\bze + \vv{v}{} + \vv{v}{}$ for some  $\vv{v}{}\in V_+$ and then $N_{\bx,\bi }$ is a monomer, for which
\[
N_{\bx,\bze + \vv{v}{} + \vv{v}{}}(\mbox{-}\vv{v}{})=\bx (\bze + \vv{v}{})=N(\vv{v}{}),
\]
\emph{i.e.}\ $N_{\bx,\bze + \vv{v}{} + \vv{v}{}}= M_{\mbox{-}\vv{v}{}:N(\vv{v}{})}$, so
\begin{equation}
\label{eqNabla}
\sum_{\{\bi \mid\; |P(\bi )\cap P(\bze )|=1\}} f(N_{\bx,\bi }) = \sum_{\vv{v}{}\in V_+}
f\left( M_{\mbox{-}\vv{v}{}:N(\vv{v}{})} \right)
 = \sum_{\vv{v}{}\in V_+}
f\left( M_{\vv{v}{}:N(\mbox{-}\vv{v}{})} \right)\,,
\end{equation}
as $V_+=-V_+$.
Combining Eqs. (\ref{eqSigmaF}), (\ref{eqDot}) and (\ref{eqNabla}), and recalling that $N_{\bx,\bze }=N$, we obtain
\begin{equation}
\begin{aligned}
\sigma (F(\bx )) =& f(N)+\sum_{(\vv{u}{},\vv{w}{})\in \CP}\left[  f\left( \dimer
{\vv{u}{}:N(\mbox{-}\vv{w}{})}{\vv{w}{}:N(\mbox{-}\vv{u}{})}\right) - f\left( \dimer{\vv{u}{}:N(\vv{u}{})}
{\vv{w}{}:N(\vv{w}{})}\right)\right]  \\
&+ \sum_{(\vv{u}{},\vv{w}{})\in \CP}
\left[f^{E}\left( \dimer{\vv{u}{}:N(\vv{u}{})}
{\vv{w}{}:N(\vv{w}{})}\right)+
f^{E}\left( \dimer{\mbox{-}\vv{w}{}:N(\vv{u}{})}
{\mbox{-}\vv{u}{}:N(\vv{w}{})}\right)\right]
+
\sum_{\vv{v}{}\in V_+}
f\left( M_{\vv{v}{}:N(\mbox{-}\vv{v}{})} \right)\,.
\end{aligned}
\end{equation}
Since $\sigma (\bx )=\sigma (F(\bx ))$, we have
\begin{equation}
\begin{aligned}
f(N) =& N(\vv{\eta}{})+
\sum_{(\vv{u}{},\vv{w}{})\in \CP}\left[ f\left( \dimer{\vv{u}{}:N(\vv{u}{})}
{\vv{w}{}:N(\vv{w}{})}\right) -  f\left( \dimer
{\vv{u}{}:N(\mbox{-}\vv{w}{})}{\vv{w}{}:N(\mbox{-}\vv{u}{})}\right)\right] +\sum_{\vv{v}{}\in V \setminus\{\vv{\eta}{}\}}f^{E}\left( H_{N(\vv{v}{})}\right)\\
&-\sum_{(\vv{u}{},\vv{w}{})\in \CP}
\left[f^{E}\left( \dimer{\vv{u}{}:N(\vv{u}{})}
{\vv{w}{}:N(\vv{w}{})}\right)+
f^{E}\left( \dimer{\mbox{-}\vv{w}{}:N(\vv{u}{})}
{\mbox{-}\vv{u}{}:N(\vv{w}{})}\right)\right]
-\sum_{\vv{v}{}\in V_+}
f\left( M_{\vv{v}{}:N(\mbox{-}\vv{v}{})} \right)\,,
\end{aligned}
\end{equation}
which concludes the proof.
\hfill\ensuremath{\square}

As an illustration, we use Theorem~\ref{maind} for $d=3$. For example, if we choose $\vv{\eta}{}=\vv{v}{1}$ and 
\[
\CP =\{ 
\lbrace \vv{v}{2},\mbox{-}\vv{v}{1}\rbrace,
\lbrace \vv{0}{},\mbox{-}\vv{v}{1}\rbrace,
\lbrace \mbox{-}\vv{v}{2},\mbox{-}\vv{v}{1}\lbrace,
\lbrace \vv{v}{3},\mbox{-}\vv{v}{1}\rbrace,
\lbrace \mbox{-}\vv{v}{3},\mbox{-}\vv{v}{1}\lbrace,
\lbrace \vv{v}{2},\mbox{-}\vv{v}{3}\rbrace,
\lbrace \vv{0}{},\mbox{-}\vv{v}{3}\rbrace,
\lbrace \vv{v}{2},\vv{v}{3}\rbrace,
\lbrace \vv{0}{},\mbox{-}\vv{v}{2}\rbrace
\}\,,
\]
then we obtain one of the $7\cdot 2^9$ formulations of the necessary and sufficient conditions for a local rule to be number-conserving.

\begin{thm}\label{main3}
Consider dimension $d=3$. A local rule $f$ is number-conserving if and only if 
for any $q_1, q_2, q_3, q_4, q_5, q_6, q_7\in Q$, it holds that
\[
\begin{aligned}
    \fff{q_1}{q_2}{q_3}{q_4}{q_5}{q_6}{q_7} = {} & q_1 +
    \fff{0}{0}{0}{q_4}{q_5}{0}{0} -\fff{0}{0}{0}{q_1}{q_2}{0}{0}+
    \fff{0}{0}{q_3}{0}{q_5}{0}{0} - \fff{0}{0}{q_1}{0}{q_3}{0}{0}  +\fff{0}{q_2}{0}{0}{q_5}{0}{0} -\fff{0}{q_1}{0}{0}{q_4}{0}{0} \\ 
&+ \fff{0}{0}{0}{0}{q_5}{q_6}{0} - \fff{0}{0}{0}{0}{q_7}{q_1}{0}  
+     \fff{0}{0}{0}{0}{q_5}{0}{q_7} -\fff{0}{0}{0}{0}{q_6}{0}{q_1}+      \fff{0}{0}{0}{q_4}{0}{0}{q_7} - \fff{0}{0}{0}{q_6}{0}{0}{q_2}  \\
     &+ \fff{0}{0}{q_3}{0}{0}{0}{q_7} -\fff{0}{0}{q_6}{0}{0}{0}{q_3} +
     \fff{0}{0}{0}{q_4}{0}{q_6}{0} - \fff{0}{0}{0}{q_7}{0}{q_2}{0} 
     +  \fff{0}{0}{q_3}{q_4}{0}{0}{0} - \fff{0}{0}{q_2}{q_3}{0}{0}{0} \\
     &+ 
   \fff{0}{0}{q_2}{0}{0}{0}{0}+\fff{0}{0}{q_6}{0}{0}{0}{0}- 2\fff{0}{0}{q_3}{0}{0}{0}{0}  +  \fff{0}{0}{0}{q_3}{0}{0}{0}  
+ \fff{0}{0}{0}{q_6}{0}{0}{0} \\
&+\fff{0}{0}{0}{q_7}{0}{0}{0}-3 \fff{0}{0}{0}{q_4}{0}{0}{0} + \fff{0}{0}{0}{0}{q_2}{0}{0} + \fff{0}{0}{0}{0}{q_3}{0}{0}+
     \fff{0}{0}{0}{0}{q_4}{0}{0} \\ 
     &+\fff{0}{0}{0}{0}{q_6}{0}{0}+
     \fff{0}{0}{0}{0}{q_7}{0}{0}-\fff{0}{0}{0}{0}{q_1}{0}{0} 
     -4\fff{0}{0}{0}{0}{q_5}{0}{0}  +   \fff{0}{0}{0}{0}{0}{q_2}{0} \\
     &- \fff{0}{0}{0}{0}{0}{q_6}{0}+
     \fff{0}{0}{0}{0}{0}{0}{q_2}+\fff{0}{0}{0}{0}{0}{0}{q_3}-
     2\fff{0}{0}{0}{0}{0}{0}{q_7}\,.
\end{aligned}
\]
%\end{equation}

\end{thm}

Although this formula is rather lengthy, it is much simpler than conditions presented by other authors (cf. \cite{Durand2003} and \cite{TI}). In three dimensions, there are $13$ number-conserving rules if $Q=\{0,1\}$. It means they are extensions of one-dimensional elementary rules equivalent to the identity rule (ECA 204 according to the classical enumeration in \cite{Wolfram83}), the shift rule (ECA 170), or the traffic rule (ECA 184), acting in each of the six possible directions.
% * <bernard.debaets@ugent.be> 2017-04-29T16:04:30.006Z:
%
% > It means they are extensions of one-dimensional elementary rules equivalent to the identity rule (ECA 204 according to the classical enumeration in \cite{Wolfram83}), the shift rule (ECA 170), or the traffic rule (ECA 184), acting in each of the six possible directions.
%
% ^.

\section{Conclusions}
In this work, we have studied $d$-dimensional CAs with the von Neumann neighborhood.
We have presented necessary and sufficient conditions for the local rule of such CA to be number-conserving. 
The greatest advantage of these conditions is that their form not only allows to decide whether a given rule is number-conserving, but also to enumerate all number-conserving rules for not too large $d$ and small $|Q|$. This is due to the fact that we do not have to
consider every rule from a huge space of $\displaystyle |Q|^{|Q|^{2d+1}}$ possible rules one by one, but only need to set at most $ (2d+1)\cdot (|Q| - 1)$ monomers and $d^2(|Q|-1)^2$ dimers and then verify whether the values given by Eq.~(\ref{eq:thm43}) belong to $Q$ for any $q_1, q_2, q_3, q_4, q_5\in Q$.
 
This method is especially useful to describe number-conserving rules satisfying some additional conditions -- in particular in cases where these additional conditions result in dependencies between monomers or dimers. This happens, for example, in the case of local rules with some kind of symmetry, the most natural one being rotation symmetry. Rotation-symmetric number-conserving  CAs (RNCAs) were studied in~\cite{TI} and~\cite{Imai2015}, where it was shown that there is no nontrivial CA that is number-conserving when $|Q|\leq 4$. Furthermore, all RNCAs with 5 states were described. These results, obtained with a considerable effort, now are a simple consequence of Eq.~(\ref{dc}). Moreover, using Eq.~(\ref{dc}), we can find all rotation-symmetric number-conserving CAs with six or even seven states. 
% * <bernard.debaets@ugent.be> 2017-04-29T16:11:23.831Z:
%
% >  which means, we do not care from which side we look at the behavior of a given CA
%
% ^.
Perhaps, one may find strongly universal automata among the newly described rules (all RNCAs with 5 states are not strongly universal).

Another type of local rules that are very interesting are passive rules, \emph{i.e.}\ rules that satisfy the condition that if a cell is surrounded by zero state cells, then its state does not change. In other words,  cells are \emph{passive} and change their value only when at least one of the neighbors has a non-zero state. Since this condition can be easily described in terms of constraints on monomers, it is possible to design a method to construct all the passive number-conserving rules of two-dimensional CAs with the von Neumann neighborhood. 

We plan to describe the solutions to the above problems in detail in a forthcoming paper.

% ref 1: delete Jul.
% ref 10: Euler-Lagrange -> capital L missing
%\printbibliography
\bibliographystyle{ieeetr}
\bibliography{numberconservingCA}

\end{document}